\documentclass{article}%
\usepackage{amsmath}%
\usepackage{amsfonts}%
\usepackage{amssymb}%
\usepackage{graphicx}

\begin{document}
\title{On the Basics of the Nonlinear Diffusion}
\author{Henrik Stenlund\thanks{The author is grateful to Visilab Signal Technologies for supporting this work.}}
\date{March 12th, 2020}
\maketitle
\begin{abstract}
This study handles spatial three-dimensional solution of the nonlinear diffusion equation without particular initial conditions. The functional behavior of the equation and the concentration have been studied in new ways. An auxiliary function for diffusion is given having an interesting relationship with the concentration. A set of new integro-differential equations is given for diffusion. 
\footnote{Visilab Report \#2020-03}
\end{abstract}
\subsection{Keywords}
diffusion, nonlinear diffusion, concentration dependence
\subsection{Mathematical Classification}
Mathematics Subject Classification: 34A34,34B15
\section{Introduction}
Diffusion is a basic natural phenomenon occurring everywhere with practically all chemical compositions thinkable, in gaseous, liquid and solid form. It happens as a self-diffusion and in mixtures. The diffusion equation is well known and has a wealth of solutions with special initial conditions, especially in one-dimensional models \cite{Churchill1972}. Diffusion is closely connected to thermal conduction problems as the equations are the same. Usually one is talking of initial value problems, not boundary value problems, as is done in this paper too. However, no specific initial value problems are treated here. In three dimensions the diffusion equation becomes awkward to solve and often numerical work is required. If a non-linearity is added things become really ugly. There are many cases where an analytical solution in closed form would be welcome for further analysis. In this paper it is attempted to give new formulas for tackling nonlinear diffusion in three dimensions. Also other observations are made of the diffusion equation. An introductory discussion is in Shewmon \cite{Shewmon1963}. The available work made by pioneers, like Boltzmann \cite{Boltzmann1894}, Matano \cite{Matano1933} and Crank \cite{Crank1979} are classic. The author has earlier presented methods for solving non-linearity based on experimental data \cite{Stenlund2011}. Any formal proofs are left out to enhance readability.
\section{The General Solution of the Diffusion Equation}
The analysis is started by the three-dimensional non-linear diffusion equation with a dependence on concentration. 
\begin{equation}
D=D(c) \label{eqn2110}
\end{equation}
Here $c$ is the concentration and $D$ is the diffusion coefficient and its functional form is supposed to be known. The diffusion equation in three dimensions will be
\begin{equation}
\frac{\partial{c}}{\partial{t}}=\nabla\cdot{(D(c)\nabla{c})} \label{eqn2120}
\end{equation}
where $t, \bar{r}$ are the time and spatial coordinate vector. No sources nor sinks are present as they will significantly complicate the analysis. One attempts to offer a solution to the initial value problem presented in the following.The analysis goes along the line of first transforming this equation to a nonlinear integral equation. Then the solution is expressed in terms of $c(\bar{r},t=0)$ and its spatial derivatives 
\begin{equation}
\frac{\partial^n{c(\bar{r},t)}}{\partial{t^n}} \label{eqn2130}
\end{equation}
which are assumed to be known. The time derivatives are formally needed but as will become clear in the following, they are actually not required, but solved instead. At $t=0$ the concentration may have either a piecewise continuous behavior or it may even be discontinuous. However, at times $t > 0$, it behaves continuously and so does its derivatives. By introducing the simple transformation
\begin{equation}
F(c)={\int_{c_0}^{c}}{D(s){{ds}}}  \label{eqn2140}
\end{equation}
with $c_0$ a constant, one will have
\begin{equation}
\nabla\cdot{(D(c)\nabla{c})}=\nabla^2{F(c)} \label{eqn2150}
\end{equation}
$F(c)$ is a scalar function of concentration and is implicitly four dimensional. Thus Eq.(\ref{eqn2120}) becomes
\begin{equation}
\frac{\partial{c}}{\partial{t}}=\nabla^2{F(c)} \label{eqn2160}
\end{equation}
 One can continue from (\ref{eqn2160}) integrating it once in terms of time
\begin{equation}
c(\bar{r},t)=c(\bar{r},t=0)+{\nabla^2}\int_0^t{F(c(\bar{r},t')dt'}  \label{eqn2200}
\end{equation}
This is likely the simplest integro-differential equation derivable for (\ref{eqn2120}). A MacLaurin power series in $t$ can be used for solving this equation. 
\begin{equation}
F(c)=\sum_{n=0}{\frac{t^n}{n!}(\frac{\partial^n{F(c(\bar{r},t))}}{\partial{t^n}})_{t=0}}  \label{eqn2220}
\end{equation}
The integral can be transformed to
\begin{equation}
c(\bar{r},t)=c(\bar{r},t=0)+\sum_{n=0}{\frac{t^{n+1}}{(n+1)!}{\nabla^2}(\frac{\partial^n{F(c(\bar{r},t))}}{\partial{t^n}})_{t=0}}  \label{eqn2230}
\end{equation}
It is important to note that this is an initial value problem. To solve general time-varying problems other methods are needed. The functions 
\begin{equation}
(\frac{\partial^n{F(c(\bar{r},t))}}{\partial{t^n}})_{t=0}
\end{equation}
can be solved from the initial differential equation (\ref{eqn2160}) at $t=0$. The first derivative becomes
\begin{equation}
\frac{\partial{F(c(\bar{r},t))}}{\partial{t}}=F'(c)\frac{\partial{c}}{\partial{t}}=F'(c){\nabla^2}F(c)
\end{equation}
The second derivative becomes
\begin{equation}
\frac{\partial^2{F(c(\bar{r},t))}}{\partial{t^2}}=F''(c)(\nabla^2{F(c)})^2+F'(c){\nabla^2}(F'(c){\nabla^2}{F(c))}
\end{equation}
Continuing in the same way for higher derivatives one can replace 
\begin{equation}
\frac{\partial{c}}{\partial{t}}
\end{equation}
by (\ref{eqn2120}). All derivatives are taken at $t=0$ and since one knows $c(\bar{r},t=0)$ and all derivatives with $\nabla$ of it and $F(c)$, one knows all the terms
\begin{equation}
(\frac{\partial^n{F(c(\bar{r},t))}}{\partial{t^n}})_{t=0}
\end{equation}
Thus (\ref{eqn2230}) is the general solution to (\ref{eqn2120}) and (\ref{eqn2200}). This result serves both as a starting point for further analysis and for numerical work in practical initial value problems. The ease of solution is dependent on the complexity of the function $F(c)$. 
\subsection{Broken Series}
In complicated cases one may approximate by breaking the series at index $N$. To estimate the resulting error, one can use Lagrange's expression for the remainder term of the series above broken at $N$.
\begin{equation}
R_N=\frac{t^{N+1}}{(N+1)!}{\nabla^2}(\frac{\partial^{N+1}{F(c(\bar{r},t))}}{\partial{t^{N+1}}})_{t=0}
\end{equation}
\section{Method for Solving the Poisson-type Differential Equation}
The following property of the three-dimensional Dirac delta function is well known
\begin{equation}
\nabla^2(\frac{1}{|\bar{r}-\bar{r}_1|})=-4{\pi}{\delta{|\bar{r}-\bar{r}_1|}} \label{eqn3000}
\end{equation}
Solving the Poisson equation
\begin{equation}
\nabla^2{V(\bar{r})}=-k(\bar{r}) \label{eqn3010}
\end{equation}
as
\begin{equation}
V(\bar{r})=\frac{1}{4{\pi}}\int\frac{k(\bar{r'})}{|{\bar{r}-\bar{r'}}|}{d\bar{r'}}  \label{eqn3020}
\end{equation}
can be done with the aid of the Dirac delta function. This can be verified by applying the Laplacian to it. One might argue that there is an additional function $\phi(\bar{r})$ involved as follows.
\begin{equation}
V(\bar{r})=\frac{1}{4{\pi}}\int\frac{k(\bar{r'})}{|{\bar{r}-\bar{r'}}|}{d\bar{r'}}+\phi(\bar{r})  \label{eqn3040}
\end{equation}
with the property
\begin{equation}
\nabla^2{\phi(\bar{r})}=0 \label{eqn3050}
\end{equation}
However, it needs to comply with the original equation simultaneously
\begin{equation}
\nabla^2{\phi(\bar{r})}=-k(\bar{r}) \label{eqn3060}
\end{equation}
Since $k()$ is arbitrary the only possibility left is
\begin{equation}
\phi(\bar{r})=0 \label{eqn3070}
\end{equation}
As an example of application one can transform equation (\ref{eqn2160}) to an integral
\begin{equation}
F(c)=-\frac{1}{4{\pi}}\int\frac{\frac{\partial{c(\bar{r}_2,t)}}{\partial{t}}}{|{\bar{r}-\bar{r}_2}|}{d\bar{r}_2}  \label{eqn3090}
\end{equation}
This method is used in many instances here.
\section{Integro-Differential Equation for the Diffusion}
The result from (\ref{eqn2120}) and opening the diffusion equation to
\begin{equation}
\frac{\partial{c}}{\partial{t}}=D'(c)(\nabla{c})^2+D(c)\nabla^2{c} \label{eqn4700}
\end{equation}
and rearranging it to get
\begin{equation}
\frac{\frac{\partial{c}}{\partial{t}}}{D(c)}-\frac{D'(c)}{D(c)}(\nabla{c})^2=\nabla^2{c} \label{eqn4710}
\end{equation}
By using the method above one will get
\begin{equation}
c(\bar{r},t)=-\frac{1}{4{\pi}}\int{\frac{d\bar{r'}}{|\bar{r}-\bar{r'}|}\big{[}\frac{1}{D(c(\bar{r'},t))}\frac{\partial{c(\bar{r'},t)}}{\partial{t}}-\nabla'{(ln(D(c(\bar{r'},t)))}\cdot{\nabla'{c(\bar{r'},t)}}\big{]}}  \label{eqn4720}
\end{equation}
This is the general integro-differential equation for nonlinear diffusion. $\nabla'$ is affecting on the $\bar{r'}$ variable only.
\section{Auxiliary Functions for the Nonlinear Diffusion Equation}
In the following the diffusion coefficient $D$ is dependent on $c$. By differentiation one has
\begin{equation}
\frac{\partial{F(c)}}{\partial{t}}=D(c)\frac{\partial{c}}{\partial{t}} \label{eqn5000}
\end{equation}
and by substituting the diffusion equation to it the result will be
\begin{equation}
\frac{\partial{F(c)}}{\partial{t}}=D(c){\nabla^2{F(c)}} \label{eqn5020}
\end{equation}
This is a nonlinear diffusion equation for the $F(c)$. It is actually an amazing equation since it is a differential equation for the integral of the diffusion coefficient in terms of the concentration variable. The $F(c)$ can be expressed as
\begin{equation}
F(c)=-\frac{1}{4{\pi}}\int\frac{d\bar{r'}}{D(c(\bar{r'},t))|\bar{r}-\bar{r'}|}{\frac{\partial{F(c(\bar{r'},t))}}{\partial{t}}}  \label{eqn5100}
\end{equation}
This is an integro-differential equation for the $F(c)$. One can make an assumption for the existence of an auxiliary function $\phi$
\begin{equation}
\phi={\nabla^2{F(c)}}  \label{eqn5110}
\end{equation}
and one can see from the above that
\begin{equation}
\phi=\frac{1}{D(c)}\frac{\partial{F(c(\bar{r},t))}}{\partial{t}}  \label{eqn5120}
\end{equation}
\begin{equation}
F=-\frac{1}{4{\pi}}\int\frac{d\bar{r'}\phi(\bar{r'},t)}{|\bar{r}-\bar{r'}|}  \label{eqn5140}
\end{equation}
Substitution to equation (\ref{eqn5000}) will give
\begin{equation}
D(c)\phi=-\frac{1}{4{\pi}}\frac{\partial}{\partial{t}}\int\frac{d\bar{r'}\phi(\bar{r'},t)}{|\bar{r}-\bar{r'}|}  \label{eqn5160}
\end{equation}
and application of a Laplacian to this produces finally
\begin{equation}
\frac{\partial{\phi}}{\partial{t}}={\nabla^2{(D(c)\phi})} \label{eqn5200}
\end{equation}
This is a non-linear diffusion equation for the auxiliary function. 
\section{Differential Equations for the Diffusion Coefficient}
It is tempting to see if the same procedures can be applied as in the preceding section to the diffusion coefficient itself.
\begin{equation}
\frac{\partial{D(c)}}{\partial{t}}=D'(c)\frac{\partial{c(\bar{r},t)}}{\partial{t}}  \label{eqn5600}
\end{equation}
\begin{equation}
\nabla{D(c)}=D'(c){\nabla{c}} \label{eqn5620}
\end{equation}
One can solve for the gradient
\begin{equation}
\nabla{c}=\frac{\nabla{D(c)}}{D'(c)} \label{eqn5640}
\end{equation}
and use the original diffusion equation to get
\begin{equation}
\frac{\partial{D(c)}}{\partial{t}}=D'(c)\nabla\cdot{(D(c)\frac{\nabla{D(c)}}{D'(c)})}  \label{eqn5680}
\end{equation}
One has obtained a differential equation for the diffusion coefficient itself, as long as the diffusion coefficient is non-linear. It is equivalent to the non-linear diffusion equation as  is very easy to see by starting to execute the differential operators. One is able to extend this thinking without a proof further to any derivative of the coefficient to
\begin{equation}
\frac{\partial{D^{(N)}(c)}}{\partial{t}}=D^{(N+1)}(c)\nabla\cdot{(D(c)\frac{\nabla{D^{(N)}(c)}}{D^{(N+1)}(c)})}  \label{eqn6690}
\end{equation}
\section{Conclusions}
The author has not seen the general three-dimensional solution (\ref{eqn2230}) given here in existing publications. Therefore, this expression is believed to be new as most of the other expressions presented. It will serve as a starting point both for analytical investigations and numerical work with various initial conditions.
 
The last part shows results for $F(c)$, the diffusion coefficient integrated in terms of the concentration having its own diffusion equation and an auxiliary function. The diffusion coefficient itself seems also to have its own differential equation. 


\begin{thebibliography}{10} 
\bibitem{Boltzmann1894}\textsc {Boltzmann, Ludwig}: \textit{Zur Integration der Diffusionsgleichung bei variabeln Diffusionscoefficienten}, Annalen der Physik, \textbf{53, 960 (1894)}
\bibitem{Matano1933}\textsc {Matano, Chujiro}: \textit{On the Relation between the Diffusion-Coefficients and Concentrations of Solid Metals (The Nickel-Copper System)}, Japanese Journal of Physics, \textbf{8, 109 (1933)}
\bibitem{Crank1979}\textsc{Crank, John}: \textit{Mathematics of Diffusion}, Clarendon Press\textbf{(1979)}
\bibitem{Shewmon1963}\textsc{Shewmon, Paul G.}: \textit{Diffusion in Solids}, McGraw-Hill New York \textbf{(1963)}
\bibitem{Churchill1972}\textsc{Churchill, Ruel}: \textit{Operational Mathematics}, McGraw-Hill Kogakusha, 3rd edition \textbf{ (1972)}, Tokyo
\bibitem{Stenlund2011}\textsc{Stenlund, H.}: \textit{Three Methods of Solution for the Concentration Dependence of the Diffusion Coefficient}, Visilab Technical Report No.2004-03. Revision 4, 05-01-2011. First published in June 2004, (www.visilab.fi/nonlinear\_diffusion.pdf) DOI: 10.13140/RG.2.1.5147.4006.
\end{thebibliography}
\end{document}